\documentclass[12pt]{amsart} 
\usepackage{multicol} 
\usepackage{amscd}

\newcommand{\Con}{\operatorname{Con}}

\newcommand{\Inv} {\operatorname{Inv}}
\newcommand{\Int}{\operatorname{Int}}
\newcommand{\cl}{\operatorname{cl}}

\newcommand{\bs}{\setminus}
\newcommand{\Ob}{\operatorname{Ob}}
\newcommand{\Mor}{\operatorname{Mor}} 
\newcommand{\M}{\mathcal{M}} \newcommand{\N}{\mathcal{N}}
\newcommand{\A}{\mathcal{A}} \newcommand{\p}{\mathcal{P}}
\newcommand{\K}{\mathcal{K}} 
 
 \newcommand{\I}{\mathcal{I}}
\newcommand{\real}{\mathbb R} \newcommand{\R}{\mathbb R}
\newcommand{\Z}{\mathbb Z} \newcommand{\ep}{\varepsilon}

\theoremstyle{plain} 
\begingroup 

\newtheorem{thm}{Theorem}[section] 
\newtheorem{cor}[thm]{Corollary}
\newtheorem{lem}[thm]{Lemma} 
\newtheorem{prop}[thm]{Proposition}
\endgroup

\theoremstyle{definition}
\newtheorem{defn}[thm]{Definition} 
\theoremstyle{remark}

\numberwithin{equation}{section}

\begin{document}

\title{Shift Equivalence and the Conley index} 
\author{John Franks and David Richeson}

\begin{abstract}
In this paper we introduce filtration pairs for isolated invariant
sets of continuous maps.  We prove the existence of filtration pairs
and show that, up to shift equivalence, the induced map on the
corresponding pointed space is an invariant of the isolated invariant
set.  Moreover, the maps defining the shift equivalence can be chosen
canonically.  Lastly, we define partially ordered Morse
decompositions and prove the existence of Morse set filtrations for
such decompositions.
\end{abstract}

\maketitle

\section{Introduction}
One would like to study the topological properties of an isolated
invariant set $S$ for a continuous map $f:X\to X$.  If $f$ is, for
instance, an Axiom A diffeomorphism with the no-cycle property and $S$
is a basic set then one may isolate $S$ between levels of a filtration
and look at the relative homology of these sets.  In general, $S$ does
not arise this way.  The discrete Conley index is a tool used to study
precisely such sets.  The ideas of Conley \cite{C} were developed by
many people using, for instance, shape theory \cite{RS}, algebraic
topology \cite{Mr} and category theory \cite{Sz}.  Due to the daunting
nature of these subjects many casual readers find Conley index theory
inaccessible.  In this paper we try to develop Conley index theory in
an intuitive fashion relying on no advanced mathematics.

We begin by introducing the filtration pair.  This topological pair
has many of the benefits of an actual filtration while possessing the
versatility of an index pair (the standard object of study in Conley
index theory).  To each filtration pair $P=(N,L)$ we may associate a
pointed topological space and an induced map $f_P:N_L\to N_L$.  This
space itself is not independent of the choice of filtration pair, but
its shift equivalence class is.  Shift equivalence is a natural
equivalence relation for non-invertible dynamical systems which is
related to topological conjugacy, but is less rigid.  It was
introduced by R. Williams for (among other things) the study of low
dimensional hyperbolic invariant sets (see \cite{W} for more details
and references).  We call the shift equivalence class, in the homotopy
category, of the map on the pointed space the {\it discrete homotopy
Conley index} of $S$.

This construction has several nice properties, not the least of which
is its relative simplicity.  We see that the basepoint is a
super-attractor; that is, $f^{-1}_P(*)$ is a neighborhood of the
basepoint $*$.  Also, unlike index pairs, filtration pairs are robust
under small $C^0$ perturbations of $f$.  If $X$ is a manifold then
$(N,L)$ may be chosen to be homeomorphic to a finite simplicial pair.
Thus, the (co)homology of the filtration pair is necessarily finite
dimensional.  Lastly, having the homotopy index one may apply various
functors to it such as homology, cohomology, $\pi_{n}$, the Leray
functor, direct limits, etc.  thus obtaining the Conley indices found
in the literature.

Next, we investigate decompositions of isolated invariant sets.  We
define attractors and repellers as a motivation so that we may define
partially ordered Morse decompositions.  We show that Morse set
filtrations exist for such decompositions.

Lastly, we discuss Szymczak's categorical definition of the Conley
index.  We show that, in fact, our construction contains exactly the
same information as his construction.

\section{Isolated Invariant Sets}

Let $U$ be an open subset of a locally compact metric space $X$ and
suppose $f: U \to X$ is a continuous map.

\begin{defn}
For any set $N \subset U$ define $\Inv^{m} N$ to be the set of $x \in
N$ such that there exists an orbit segment $\{x_{n}\}_{-m}^{m} \subset
N$ with $x_{0}=x$ and $f(x_{n})=x_{n+1}$ for $n=-m,\ldots, m-1.$ We
will call a complete orbit containing $x$ a {\em solution through x.}
More precisely, if $\sigma: \Z \to N$ is given by $\sigma(n) = x_n$
and $x_{0}=x$ and $f(x_{n})=x_{n+1}$ for all $n$, we will call
$\sigma$ a solution through $x.$ We define the {\em maximal invariant
subset,} $\Inv N$ to be $\Inv^\infty N$, the set of $x \in N$ such
that there exists a solution $\sigma$ with
$\{\sigma(n)\}_{-\infty}^{\infty} \subset N$ and with $\sigma(0)=x$.
\end{defn}

Note that {from} the definition it is clear that $f(\Inv N) = \Inv N.$
The following basic property of $\Inv N$ is trivial if $f$ is
one-to-one, but requires proof in the case of a general continuous
map.

\begin{prop}\label{Invm}
If $U$ is an open subset of a locally compact metric space $X, \ f: U
\to X$ is continuous and $N$ is a compact subset of $U$, then
\[
\Inv N = \bigcap_{m=0}^\infty \Inv^m N
\]
\end{prop}

\begin{proof}
Trivially $\Inv N \subset \bigcap_{m=0}^\infty \Inv^m N$ so we need
only prove the reverse inclusion.  If $x \in \bigcap_{m=0}^\infty
\Inv^m N$ then we must show that $x \in \Inv N.$

Let $X_1 = f^{-1}(x) \cap N$ and, for $k>1$, inductively define $X_k =
f^{-1}(X_{k-1}) \cap N$ so each $X_k$ is non-empty, compact and
$f(X_k) \subset X_{k-1}.$ Then define
\[
Y_k = \bigcap_{n \ge 1} f^n( X_{n+k})
\]
Since $\{f^n(X_{n+k})\}$ is a nested sequence of non-empty compact
subsets of $X_k$, it follows that $Y_k$ is a non-empty compact subset
of $X_k$.  So $f(Y_1) = \{x\}$ and for $k>1$,
\[
f(Y_k) = \bigcap_{n \ge 1} f^{n+1}( X_{n+k}) = \bigcap_{n \ge 2} f^n(
X_{n+k-1}) = Y_{k-1}
\]
Thus every point of $Y_{k-1}$ has an $f$ pre-image in $Y_k$.  Define
$x_{-1}$ to be any point of $Y_1$ (so $f(x_{-1}) = x$) and inductively
define $x_{-k}$ to be a point in $Y_k$ such that $f(x_{-k}) =
x_{-k+1}.$ If we also define $x_k = f^k(x)$ for $k \ge 0$ then
$\sigma: \Z \to N$, given by $\sigma(n) = x_n$, is a solution through
$x.$ It follows that $x \in \Inv N.$
\end{proof}

\begin{defn}
A compact set $N$ is called an {\em isolating neighborhood}\footnote{This definition of isolating neighborhood differs slightly
{from} that used by R. Easton in \cite{E}, for example.  There, $N$ is
called an isolating neighborhood provided $\cap_{m \in \Z} f^m(N)$ is
in the interior of $N$.  For a map which is not one-to-one the sets
$\Inv N$ and $\cap_{m \in \Z} f^m(N)$ may not coincide, as is pointed
out in \cite{E}.} if $\Inv N\subset \Int N$.  A set $S$ is called an
{\em isolated invariant set} if there exists an isolating neighborhood
$N$ with $S=\Inv N$.  A compact set $N$ is called an {\em isolating
block} if $f(N)\cap N\cap f^{-1}(N)\subset\Int N$.
\end{defn}

Notice that every isolating block $N$ is an isolating neighborhood for
the set $S = \Inv N$.  The converse is not true.

Suppose $N$ is an isolating neighborhood.
\begin{defn}
We define the {\em exit set} of $N$ to be \[ N^{-}:=\{x\in N: f(x)
\notin \Int N\}.
\]
\end{defn}

That is, $N^-$ is the subset consisting of those points in $N$ whose
image is in the complement of the interior of $N$.

\section{Filtration Pairs}

Ideally, for the study of an isolated invariant set $S$ we would like
a pair $(N,L)$ of spaces forming part of a filtration; i.e.  we would
like to have $L\subset N$ compact, $S = \Inv (N \bs L), \ f(L) \subset
\Int L,$ and $f(N) \subset \Int N$.  This is much too strong, however,
as we often need to study $f$ on a small neighborhood of $S$ and there
may not exist a forward invariant small neighborhood $N$ of $S$.

To motivate the proper modification of the notion of a filtration
pair, we will look at the definition of a filtration.  Consider a map
$f$ defined on a neighborhood of a pair of compact sets $(N,L)$ equal
to the closures of their interiors with $f(N)\subset N$ and
$f(L)\subset L$.  Then the condition $f(L) \subset \Int L$ is
equivalent to $f(L) \cap \cl (N \bs L) = \emptyset.$ Also, given this,
the condition $f(N) \subset \Int N$ is equivalent to $f(\cl(N \bs L))
\subset \Int N.$ We take these two conditions for the case when $N$
and $L$ are not forward invariant.  Notice that $f(\cl(N \bs L))
\subset \Int N$ is equivalent to requiring that $L$ be a neighborhood
of $N^{-}$.  This requirement also asserts that any point of $N$ which
is mapped by $f$ outside of $N$ must be in $L$.  Thus we have the
following definition.

\begin{defn}\label{filpair}
Let $S$ be an isolated invariant set and suppose $L \subset N$ is a
compact pair contained in the interior of the domain of $f$.  The pair
$(N,L)$ is called a {\em filtration pair} for $S$ provided $N$ and $L$
are each the closures of their interiors and
\begin{enumerate}
\item
$\cl (N\bs L)$ is an isolating neighborhood of $S$,
\item
$L$ is a neighborhood of $N^{-}$ in $N$, and
\item
$f(L) \cap \cl (N\bs L) = \emptyset$.
\end{enumerate}
\end{defn}

The definition of filtration pair here is close to, but somewhat more
general than, the concept of filtration pair introduced in \cite{BF}
where the intent was to study dynamical systems satisfying stronger
hypotheses than we now wish to assume.  Also, our definition of
filtration pair is close to, but somewhat less general than, the
definition of {\it index pair} introduced by Conley and studied by
several authors (see \cite{E}, \cite{Mr}, and \cite{Sz}).  It is
significantly less general than the definition of index pair found in
\cite{RS}.

We have another property that is analogous to one for a filtration.
If $(N_1, L)$ and $(N_2, L)$ are filtration pairs for isolated
invariant sets $\Inv(N_{1}\bs L)$ and $\Inv(N_{2}\bs L)$ and $N_1
\subset \Int N_2$ then $(N_2, N_1)$ is a filtration pair for
$\Inv(N_{2}\bs N_{1})$.

Suppose $f:X\to X$ is a continuous map and $x\in X$.  We define the {\em
$\omega$-limit set} of $x$ to be
\[
\omega(x)=\underset{N>0}{\bigcap}\cl\Big(\underset{n>N}{\bigcup}f^n(x)\Big
).\] We will need the following easy proposition in our study of
filtration pairs.

\begin{prop}\label{omega}
If $P = (N,L)$ is a filtration pair for $S$ and $x \in N$ satisfies
$f^n(x) \in N \bs L$ for all $n > 0$ then $\omega(x) \subset S.$
\end{prop}

\begin{proof}
If $n_i \to \infty$ and the sequence $\{f^{n_i}(x)\}$ converges to $z \in
\cl (N \bs L)$ then for any $k \in \Z$ the sequence
$\{f^{n_i+k}(x)\}$ is well defined for $i$ sufficiently large and
converges to a point $z_k \in \cl (N \bs L)$.  By continuity of $f$,
if $k > 0$ then $z_k = f^k(z)$ and if $k <0$ then $z \in f^{-k}(z_k)$.
It follows that $z \in \Inv \cl(N \bs L) = S.$
\end{proof}

Our next objective is to prove the existence of isolating blocks and
filtration pairs inside any neighborhood of $S$.

Recall that for any $\varepsilon>0$, a sequence $\{x_{n}\}^{p}_{q}$ is
called an $\varepsilon$-{\em chain} provided
\[
d(f(x_{n}),x_{n+1})<\varepsilon
\]
for all $n=q,\ldots, p-1$.

\begin{defn}
For any isolating neighborhood $N$ of an isolated invariant set $S$ and any
$\varepsilon>0$ define the {\em $\varepsilon$-chain neighborhood} of
$S$ relative to $N$ to be the set $C_{\varepsilon}(N,S)$ of all $x \in
N$ such that, for some $k$, there exists an $\varepsilon$-chain
$\{x_{n}\}_{-k}^{k}\subset N$ with $x=x_{0}$ and $x_k, x_{-k} \in S.$
\end{defn}

In the next proposition we show that the sets $C_{\varepsilon}(N,S)$
for different $\varepsilon >0$ provide a neighborhood basis for $S$.

\begin{prop}\label{ep_nbhd}
Let $S$ be an isolated invariant set with isolating neighborhood $N$.
For any neighborhood $V$ of $S$ in $N$ there is an $\varepsilon >0$
with $C_{\varepsilon}(N,S) \subset V$.
\end{prop}

\begin{proof}
We may assume $V$ is open.  {From} Proposition~\ref{Invm} it follows
that there is an $m>0$ with $\Inv^m N \subset V.$ 
Our proof is by contradiction.  Suppose there is no $\ep > 0$ with
$C_{\ep}(N,S) \subset V$.  Then 
every $\ep > 0$ there is
an $x(\ep) \in C_{\varepsilon}(N,S) \setminus V,$ i.e.  there is an
$\varepsilon$-chain $\{x_{n}(\ep)\}_{-k(\ep)}^{k(\ep)}\subset N$ with
$x(\ep)=x_0(\ep)$ and $x_{k(\ep)}(\ep), x_{-k(\ep)}(\ep) \in S.$
Truncating these when $k(\ep)>m$ or extending them to be orbit
segments if $k(\ep)<m$ we produce for each $\ep >0$ a $\ep$-chain
$\{x_{n}(\ep)\}_{-m}^{m}\subset N$ with $x_0(\ep) \notin V.$ Choosing
a sequence $\ep_i$ converging to $0$ as $i \to \infty$, we may assume
that each of the sequences $\{x_j(\ep_i)\},\ -m \le j \le m$ has a
limit, say $x_j$.  Clearly $x_j \in N$ and $x_0 \in N \setminus V.$ Also it
is clear by continuity of $f$ that $f(x_j) = x_{j+1}$ for $-m \le j
<m.$ Hence $x_0 \in \Inv^m N \subset V$ which is a contradiction.
\end{proof}

We now give a proof of a result due to R. Easton \cite{E} which
asserts the existence of isolating blocks.  We essentially reproduce
Easton's proof, because it is short, elementary and quite elegant.

\begin{prop}
Every neighborhood of an isolated invariant set $S$ contains an
isolating block.  In particular if $N$ is an isolating neighborhood
for $S$, then for every sufficiently small $\varepsilon>0, \
\cl(C_{\varepsilon}(N,S))$ is an isolating block.
\end{prop}

\begin{proof}
Let $B = \cl(C_{\varepsilon}(N,S)).$ If $B$ is not an isolating block
for $S$ then there is a point $x$ such that $x, f(x)$ and $f^2(x)$ are
all in $B$, but $f(x) \notin \Int B.$ By continuity of $f$ there
exists $\delta>0$ such that if $d(x,y)<\delta$ then
$d(f(x),f(y))<\varepsilon$.  Since $x \in B$ and $f^2(x) \in B$ there
are $\varepsilon$-chains $\{x_{n}\}_{0}^{k}$ with $x_0 \in S$ and
$d(x_k, x) < \delta$, and $\{y_n\}_{0}^{m}$ with $y_m \in S$ and
$d(y_0, f^2(x)) < \ep$.  Then the sequence $x_0, x_1, \dots, x_k, f(x),
y_0, \dots, y_m$ is an $\varepsilon$-chain {from} $S$ to $S$
containing $f(x)$.  This implies $f(x) \in C_{\varepsilon}(N,S) = \Int
B,$ which is a contradiction.
\end{proof}

\begin{thm} \label{fpexist}
If $N$ is an isolating block and $L$ is any sufficiently small compact
neighborhood of $N^-$ in $N$ then $(N,L)$ is a filtration pair.
Moreover there is a neighborhood of $f$ in the $C^0$ topology such
that for any $\tilde f$ in this neighborhood, $\tilde S = \Inv
(N\setminus L, \tilde f)$ is an isolated invariant set and $(N,L)$ is
a filtration pair for $\tilde S$.
\end{thm}

\begin{proof}

To check property (1) of Definition~\ref{filpair} it suffices to show
$S \subset N\setminus L.$ Clearly $f(N^-)$ is in the complement of $\Int N$
and hence $S \cap f(N^-) =\emptyset.$ Thus if $L$ is a sufficiently
small neighborhood of $N^-$, we have $S \cap f(L) =\emptyset.$ Since
$S$ is invariant $S \cap L$ could only be non-empty if $S \cap f(L)$
were.  Hence $S \cap L =\emptyset$ and $S \subset N\setminus L.$

Property (2) of Definition~\ref{filpair} is trivially satisfied.
To check property (3) note first that since
\[
N \setminus N^- \subset N \cap f^{-1}(\Int N) \subset N \cap f^{-1}(N),
\]
we have $\cl (N \setminus N^-) \subset N \cap f^{-1}(N).$

Also $f(N^-) \subset f(N) \cap (\Int N)^c$ where the superscript $c$
indicates complement.  It follows that
\[
f(N^-) \cap \cl (N \setminus N^-) \subset (f(N) \cap (\Int N)^c) \cap N \cap
f^{-1}(N) = \emptyset,
\]
because $f(N) \cap N \cap f^{-1}(N) \subset \Int N.$

Since $f(N^-)$ and $\cl (N \setminus N^-)$ are disjoint it is clear that if
$L$ is a sufficiently small neighborhood of $N^-$ then $f(L)$ and $\cl
(N \setminus L)$ are also disjoint.  Thus property (3) of
Definition~\ref{filpair} is satisfied.

Finally we check the fact that $(N,L)$ is a filtration pair for maps
which are $C^0$ close to $f$.  First we note that $N$ is an isolating
block for nearby maps.  If this were not the case then there would
exist a sequence $\{f_n\}$ converging to $f$ and a sequence $\{x_n\}$
with $x_n \in f_n^{-1}(N) \cap N \cap f_n(N) \cap (\Int N)^c.$
Choosing a subsequence we can assume $\{x_n\}$ converges to $z \in
f^{-1}(N) \cap N \cap f(N) \cap (\Int N)^c$ which would contradict the
fact that $N$ is an isolating block for $f$.

Similarly if $N_n^{-}=\{x\in N: f_n(x) \notin \Int N\}$ and $x_n \in
N_n^{-} \cap (\Int L)^c$ then a subsequence of $\{x_n\}$ will converge
to a point $z \in N^{-} \cap (\Int L)^c$ contradicting the hypothesis
that $L$ is a neighborhood of $N^{-}.$ Thus $(N,L)$ is a filtration
pair for any map sufficiently $C^0$ close to $f$.
\end{proof}

In \cite{RS} Robbin and Salamon give a robustness theorem for index
pairs (Theorem 5.3).  We remarked following Definition \ref{filpair}
that they had a very general definition of index pairs.  It should be
noted that the robust index pairs constructed in the proof of this
theorem are not filtration pairs nor are they index pairs in the
sense this term is used in \cite{E}, \cite{Mr}, and \cite{Sz}.

\begin{thm}\label{small_fp}
Let $U$ be an open subset of an $n$ dimensional manifold $M$ and
suppose $f: U \to M$ is a continuous map with an isolated invariant
set $S$.  Inside any neighborhood of $S$ there exists a filtration
pair $(N,L)$ such that $N$ is an $n$ dimensional manifold with
boundary and $\partial L$ is an $n-1$ dimensional submanifold of $M$.
In particular $(N,L)$ is homeomorphic to a finite simplicial pair.
\end{thm}

\begin{proof}
We may assume the given neighborhood $V$ is open in $M$ and $\cl V$ is
a compact subset of $M$.  Let $N_1$ be an isolating block for $S$
inside $V$.

Let $\varphi:\cl V\to\real^{+}$ be a smooth function which vanishes
exactly on $N_1.$ For any regular value $\delta>0$
$N=\varphi^{-1}([0,\delta])$ is a compact manifold with boundary.  If
we choose $\delta$ small enough $N$ will be an isolating block for
$S$.

Now, let $\psi:\cl(V)\to\R^{+}$ be a smooth map which vanishes exactly
on $N^{-}$.  For any regular value $\eta>0$ the space
$\psi^{-1}([0,\eta])$ is a manifold with boundary.  If $\eta>0$ is
small enough and $L_0:=\psi^{-1}([0,\eta])$, then $L_0$ is a small
neighborhood of $N^-$.  Modifying $L_0$ slightly we may assume that
the boundary of $N$ and the boundary of $L_0$ intersect transversely.
If we define $L := N \cap L_0$ and $\eta$ was chosen sufficiently
small then it follows {from} Theorem~\ref{fpexist} that $(N, L)$ is a
filtration pair.  Triangulating $N$ with $L$ as a subcomplex shows
that $(N, L)$ is homeomorphic to a finite simplicial pair.
\end{proof}

\begin{thm}\label{quotient}
Let $P=(N,L)$ be a filtration pair for $f$ and let $N_L$ denote the
quotient space $N/L$ where the collapsed set $L$ is denoted $[L]$ and
is taken as the base-point.  Then $f$ induces a continuous base-point
preserving map $f_P : N_{L}\to N_{L}$ with the property $[L]\subset
\Int f_P^{-1}([L])$.  This map will be called the pointed space map
associated to $P$.
\end{thm}

Note that in case $L$ is empty then $N/L$ is defined to be the
disjoint union of $N$ and the single point space $[L] := [\emptyset].$

\begin{proof}
Let $p:N \to N/L$ be the quotient map and define $f_{P}([L]) = [L]$
and $f_{P}(x) = p(f(x))$ otherwise, where we have identified $N_L \setminus
\{[L]\}$ with $N \setminus L.$ By Definition~\ref{filpair} $f(L)$ is
disjoint {from} $\cl (N \setminus L)$.  Hence if $K$ is a sufficiently small
neighborhood of $L$ in $N$, then $f(K)$ is disjoint {from} $\cl (N \setminus
L)$.  Hence $f_P(x) = [L]$ for all $x \in p(K).$

It is immediate that $f_P$ is continuous on $N_L \setminus \{[L]\}$ since it
is the composition of continuous functions there.  To check continuity
at $[L]$ note that if $\{x_n\}$ is a sequence in $N_L \setminus \{[L]\}$
converging to $[L]$ then $\{p^{-1}(x_n)\}$ is eventually in $K$ so
$\{f_P(x_n)\}$ is a sequence which is eventually constant and equal to
$[L].$
\end{proof}

The pointed space map $f_P$ depends on more than the invariant set
$S$.  Indeed even its homotopy type and the homotopy type of $N_L$
depend on the choice of filtration pair $P = (N,L).$ Our next section
is devoted to an important equivalence relation for which all choices
of filtration pair will give equivalent associated pointed space maps.

\section{Shift equivalence}

Suppose $\K$ is a category.  Let $X,X^{\prime}$ be objects in $\K$ and
$f:X\to X$, $g:X^{\prime}\to X^{\prime}$ be endomorphisms.  We say
that $(X,f)$ and $(X^{\prime},g)$ are {\em shift equivalent},
$f\sim_{s}g$, if there exist $m\in\Z^{+}$, $r:X\to X^{\prime}$ and
$s:X^{\prime}\to X$ such that the diagrams

\begin{multicols}{2}
\begin{center}
$\begin{CD} X @>f>> X\\ @VrVV @VrVV\\ X^{\prime} @>g>> X^{\prime}
\end{CD}$\\

$\begin{CD} X^{\prime} @>g>> X^{\prime}\\ @VsVV @VsVV\\ X @>f>> X
\end{CD}$
\end{center}
\end{multicols}
\noindent commute and $r \circ s=g^{m}$ and $s \circ r=f^{m}$.  The
integer $m$ is called the {\em lag}.

Shift equivalence is a natural and dynamically significant equivalence
relation for maps.  Note that if $f$ and $g$ are homeomorphisms (i.e.
invertible) then they are shift equivalent if and only if they are
topologically conjugate.  If they are shift equivalent a conjugacy is
given by $h = r \circ f^{-m} = g^{-m} \circ r$ and $h^{-1} = s.$

We wish to prove that for any two filtration pairs of an isolated
invariant set $S$ the corresponding maps on the pointed spaces are
shift equivalent.  We do this by first proving it for two special
cases in the following lemmas.

\begin{lem}\label{L1}
Suppose $P^{\prime} = (N,L^{\prime})$ and $P = (N \cup L, L)$ are
filtration pairs for $S$ and that $L^{\prime} \subset L$ and $f(L)
\subset \Int L.$ Then the induced maps, $f_{P^{\prime}}$ and $f_P$, on
the corresponding pointed spaces are shift equivalent.
\end{lem}

\begin{proof}

Let $Q = N \cup L.$ There is a continuous function $r: N_{L^{\prime}}
\to Q_L$ given by $r([L^{\prime}]) = [L]$ and $r(x) = p_2(x)$
otherwise, where we have identified $N_{L^{\prime}} \setminus
\{[L^{\prime}]\}$ with $N \setminus L^{\prime}$ and $p_2:Q \to Q_{L}$ is the
quotient map.  Clearly $r \circ f_P = f_{P^{\prime}} \circ r.$

Since $L^{\prime}$ contains the exit set of $N$ and $S$ is disjoint
{from} $L$, there is an $n >0$ such that for every $x \in N \cap L$ we
have $f^k(x) \in L^{\prime}$ for some $k < n.$

Define a map $s:Q_L \to N_{L^{\prime}}$ by $s([L]) = [L^{\prime}]$ and
$s(x) = f^n_{P^{\prime}} (p_1(x))$ otherwise, where we have identified
$Q_L \setminus \{[L]\}$ with $Q \setminus L = N \setminus L$ and $p_1:
N \to N_{L^{\prime}}$ is the quotient map.  Clearly $s \circ f_P =
f_{P^{\prime}} \circ s.$ 

It is immediate that $s$ is continuous on $Q_L \setminus \{[L]\}$
since it is the composition of continuous functions there.  To check
continuity at $[L]$ note that there is a neighborhood $V$ of $L$ in
$Q$ such that $f(V) \subset \Int L$ and hence for any point $z \in V$,
we have $f^k(z)\in L^{\prime}$ for some $k \le n$.  Thus for any point
$z \in N \cap V$, we have $f^n_{P^{\prime}} (p_1(z)) = [L^{\prime}].$
So $s(x)$ is constant on a neighborhood of $[L]$ in $Q_L$ and
hence continuous at $[L]$.

It is also clear {from} the definitions that $s \circ r =
f_{P^{\prime}}^n$ and $r \circ s = f_P^n.$
\end{proof}

\begin{lem}\label{L2}
Suppose $P = (N,L)$ and $P^\prime = (N^\prime, L)$ are filtration
pairs for $S$ and that $N \setminus L \subset N^\prime \setminus L$ and $f(L)
\subset \Int L.$ Then the induced maps, $f_{P}$ and $f_{P^\prime}$, on
the corresponding pointed spaces, are shift equivalent.
\end{lem}

\begin{proof}
There is a continuous function $r: N_{L} \to N^{\prime}_L$ given by
$r([L]) = [L]$ and $r(x) = p^\prime(x)$ otherwise, where we have
identified $N_{L} \setminus \{[L]\}$ with $N \setminus L$ and $p^\prime:N^\prime
\to N^{\prime}_L$ is the quotient map.  Clearly $r \circ
f_{P^\prime} = f_P \circ r.$

It follows {from} Proposition \ref{omega} that there is an $n>0$ such
that $f^n(N^\prime) \subset \Int N.$ Define $s: N^\prime_L \to N_{L}$
by $s([L]) = [L]$ and $s(x) = p(f^n(x))$ otherwise, where we have
identified $N^\prime_L \setminus \{[L]\}$ with $N^\prime \setminus L$ and $p:
N \to N_L$ is the quotient map.  Clearly $s \circ f_{P^\prime} =
f_P \circ s.$

It is immediate that $s$ is continuous on $N^\prime_L \subset \{[L]\}$
since it is the composition of continuous functions there.  To check
continuity at $[L]$ note that if $\{x_n\}$ is a sequence in
$N^\prime_L$ converging to $[L]$ then $\{s(x_n)\}$ is a sequence which
is eventually constant and equal to $[L] \in N_L$.
\end{proof}

\begin{thm}
\label{equivthm}
Suppose $P = (N,L)$ and $P^{\prime} = (N^{\prime}, L^{\prime})$ are
filtration pairs for $S$.  Then the induced maps, $f_{P}$ and
$f_{P^{\prime}}$, on the corresponding pointed spaces, are shift
equivalent.
\end{thm}

\begin{proof}
By Proposition~\ref{ep_nbhd} we may choose $\ep>0$ sufficiently small
that $C_\ep(N \setminus L, S) = C_\ep(N^{\prime} \setminus
L^{\prime},S) \subset \Int (N \setminus L)\cap \Int (N^{\prime} 
\setminus L^{\prime})$.  Denote $\cl C_\ep(N \setminus L,S)$ 
by $B$.  We will show that if $B_0$ is a sufficiently small compact
neighborhood in $B$ of the exit set $B^-$ then $P_0=(B,B_0)$ is a
filtration pair with the property that $f_{P_0}$ is shift equivalent
to both $f_{P}$ and $f_{P^{\prime}}$.  The argument is the same in
both cases so we will give the proof for $f_{P}$.

As before we will identify $N \setminus L$ with
$N_{L}\setminus\{[L]\}$ and hence we may consider $B$ as a subset of
$N_{L}$ and consider the map to be $f_{P}$.  Recall that $B^-$ is the
set of $x \in B = C_\ep(N,S)$ such that $f_{P}(x) \notin \Int B.$ Thus
{from} the definition of $C_\ep(N,S)$ it is clear that $x\in B^-$
implies there is no $\ep$-chain {from} $f_{P}(x)$ to $S$.
Consequently by Proposition~\ref{omega} the $\omega$-limit set of such
an $x$ can only be the point $[L]$.  {From} this it follows that there
is an $n>0$ such that for any $x \in B^-$ we have $f^n_{P}(x) = [L].$
If $B_0$ is a sufficiently small neighborhood of $B^-$ in $B$, the
same is true for points of $B_0.$

Define ${K} \subset N_{L}$ to be $\cl (\Int f^{-n}_{P}([L])).$ Then
$f_{P}({K}) \subset \Int {K}$ because $f^{-1}_{P}([L])$ is a
neighborhood of $[L]$ and any sequence converging to a point of
$f_{P}({K})$ will have an image under $f^{n-1}_{P}$ converging to
$[L]$.  So any such sequence will eventually be in ${K}.$ We also note
that $B_0 \subset \Int {K}$.

Thus $R = (N_{L},{K})$ and $Q = (B\cup {K},
{K})$ are filtration pairs for $S$.  We complete the proof by
showing $f_R: N_{L}/{K} \to N_{L}/{K}$ is shift
equivalent to both $f_P: N_{L} \to N_{L}$ and $f_{P_0}:
B_{B_0} \to B_{B_0}.$ By Lemma~\ref{L1} $f_{P_0}$ is shift equivalent
to $f_{Q}$ and by Lemma~\ref{L2} $f_{Q}$ is shift equivalent to
$f_{R}$.  Thus $f_{R}$ is shift equivalent to $f_{P_0}$

To see that $f_{R}$ is shift equivalent to $f_{P}$ we observe
that if $p: N \to N_{L}$ is the quotient map then we can identify the
associated pointed spaces (and maps) of the filtration pairs $(N_{L},
{K})$ and $(N, p^{-1}({K}))$.  By Lemma~\ref{L1} $f_{P}$
is shift equivalent to the associated pointed space map for the pair
$(N, p^{-1}({K}))$ which we have identified with $f_{R}$.
\end{proof}

One could ask if the shift equivalence constructed in Theorem
\ref{equivthm} is unique.  Clearly it is not.  For any two spaces and
a given lag there may be many different maps $r$ and $s$ giving a
shift equivalence.  The maps constructed in Theorem \ref{equivthm} do
have one distinguishing property.  Loosely speaking, each map takes a
point $x$ either to $f^k(x)$ or to the basepoint.  We characterize
this property as follows.

\begin{defn}
If $r: N_{L} \to N^\prime_{L^\prime}$ and $s: N^\prime_{L^\prime} \to
N_{L}$ define a shift equivalence of lag $m$, we will say that this is
a {\em standard shift equivalence} provided there exists a non-negative
integer $k\le m$ such that for every $x \in N \setminus L$, either
\begin{enumerate}
\item
$f^k(x) \in N^\prime \setminus L^\prime$, and $r(x) =
p^\prime(f^k(x))$, where $p^\prime: N^\prime \to N^\prime_{L^\prime}$
is the quotient map, or
\item
$r(x) = [L^\prime]$.
\end{enumerate}
We call $k$ the {\em semi-lag} of $r$.  Also we require that the map
$s$ has the analogous property with semi-lag $m-k$.
\end{defn}

Still, there may be many standard shift equivalences between two
pointed spaces.  But, in Proposition \ref{canon} we show that for
large enough semi-lag there exist canonical choices for each
semi-conjugacy.

\begin{lem}\label{stand_exist}
The maps $r$ and $s$ defined in Theorem \ref{equivthm} give a standard
shift equivalence.
\end{lem}

\begin{proof}
It is immediate {from} the proofs of Lemmas \ref{L1} and \ref{L2} that
this result holds in the two special cases covered by those lemmas.
But in the general case proved in Theorem \ref{equivthm} the
semi-conjugacies constructed were all compositions of semi-conjugacies
arising in these two special cases.  It is clear {from} the definition
that the composition of semi-conjugacies of standard shift
equivalences is again standard, so the result follows.
\end{proof}

\begin{lem}\label{standlem}
Suppose $r$ and $s$ define a standard shift equivalence of lag $m$
{from} $f_P: N_{L} \to N_{L}$ to $f_{P^\prime}: N^\prime_{L^\prime}
\to N^\prime_{L^\prime}$ with semi-lags $k$ and $m-k$ respectively.
Then if $r^\prime$ and $s^\prime$ define another standard shift
equivalence between $f_P$ and $f_{P^\prime}$ with the same semi-lags,
\[
r\circ f_P^m = r^\prime\circ f_P^m \text{ and } s \circ f_{P^\prime}^m
= s^\prime\circ f_{P^\prime}^m.
\]
\end{lem}

\begin{proof}
Let $x \in N_{L}$.  If $f_P^m(x) = [L]$ then clearly $r\circ f_P^m(x)
= r^\prime\circ f_P^m(x).$ On the other hand, if $f_P^m(x) \ne [L]$
then $r(x) \ne [L^\prime]$ since this would imply $f_P^m(x)= s(r(x)) =
s([L^\prime]) = [L].$ Likewise we have $r^\prime(x) \ne [L^\prime]$.
Hence it follows that $r(x) = p^\prime(f^k(x)) = r^\prime(x).$ {From}
this and the fact that $r$ and $r^\prime$ are semi-conjugacies, it is
immediate that $r(f_P^m(x)) = r^\prime(f_P^m(x)).$

The proof that $s \circ f_{P^\prime}^m = s^\prime\circ f_{P^\prime}^m$
is similar.
\end{proof}

The following result is similar to Theorem 6.3 of \cite{RS} in which an
analogous result was proved for invertible dynamical systems.

\begin{prop} \label{canon}
For any two filtration pairs $P=(N,L)$ and
$P^\prime=(N^\prime,L^\prime)$ there exists a number $M(P,P^\prime)\ge
0$ and maps $r^k_{P^\prime P}:N_L\to N^\prime_{L^\prime}$ defined for
$k\ge M(P,P^\prime)$, with the following properties
\begin{enumerate}
\item
$M(P,P^\prime)=M(P^\prime,P)$ and
\item
$r^k_{PP}=f^k_P$.
\end{enumerate}
If $P^{\prime\prime}=(N^{\prime\prime},L^{\prime\prime})$ is another
filtration pair then
\begin{enumerate}
\item[3.]
$M(P,P^{\prime\prime})\le
M(P,P^{\prime})+M(P^\prime,P^{\prime\prime})$ and
\item[4.]
$r^{k_1+k_2}_{P^{\prime\prime}P }=r^{k_1}_{
P^{\prime\prime}P^\prime}\circ r^{k_2}_{P^{\prime}P}$.
\end{enumerate}
\end{prop}
\begin{proof}
By Lemma \ref{stand_exist} we know that $f_P:N_L\to N_L$ and
$f_P^\prime:N^\prime_{L^\prime}\to N^\prime_{L^\prime}$ are shift
equivalent via a standard shift equivalence.  Let $m$ be the minimal
lag of all standard shift equivalences.  Take $M(P,P^\prime)=2m$.
Defining $M(P,P^{\prime})$ in such a way it is clear that properties
(1) and (3) hold.

Let $k_1\ge M(P,P^\prime)$.  We construct $r^{k_1} _{P^\prime P}:N_L\to
N^\prime_{L^\prime}$ as follows.  Let $r:N_L\to N^\prime_{L^\prime}$
and $s:N^\prime_{L^\prime}\to N_L$ be a standard shift equivalence
realizing the minimal lag $m$ and having semi-lags $k$ and $m-k$
respectively.  Define \[r^{k_1}_{P^\prime P}=r\circ f_P^{k_1-k}.\]

First, notice that since $k_1-k\ge m$ Lemma \ref{standlem} implies
that $r^{k_1}_{P^\prime P}$ is independent of $r$ and $s$.  Also, if
$P=P^\prime$ then $M(P,P)=0$ and taking $r=id$ we see that (2) is
satisfied.

Let $k_2\ge M(P^\prime,P^{\prime\prime})$ be given.  Let
$r^{\prime\prime}:N_{L}\to N^{\prime\prime}_{L^{\prime\prime}}$ and
$s^{\prime\prime}: N^{\prime\prime}_{L^{\prime\prime}}\to N_{L}$ form
a standard shift equivalence with semi-lags $k^{\prime\prime}$ and
$m^{\prime\prime}-k^{\prime\prime}$ respectively achieving minimal lag
$m^{\prime\prime}$.  Observe that $r^\prime\circ r$ and $s\circ
s^\prime$ form a standard shift equivalence between $f_P$ and
$f_{P^{\prime\prime}}$ with semi-lags $k+k^\prime$ and
$m+m^\prime-k-k^\prime$ respectively.  Since $k_1+k_2\ge
M(P,P^{\prime\prime})=2m^{\prime\prime}$ it follows that
$k_1+k_2-k^{\prime\prime}\ge m^{\prime\prime}$, the lag of the shift
equivalence given by $r^{\prime\prime}$ and $s^{\prime\prime}$.
Similarly, $k_1+k_2-k-k^{\prime}\ge m+m^\prime,$ the lag of
$r^\prime\circ r$ and $s\circ s^\prime$.  Thus, by Lemma
\ref{standlem} we have

\begin{align*}
r^{k_1+k_2}_{P^{\prime\prime} P}&=r^{\prime\prime}\circ
f_P^{k_1+k_2-k^{\prime\prime}}\\ &= r^\prime\circ r\circ
f_P^{k_1+k_2-k-k^\prime}\\ &= r^\prime\circ
f_{P^\prime}^{k_2-k^\prime}\circ r\circ f_P^{k_1-k}\\
&=r^{k_2}_{P^{\prime\prime} P^\prime}\circ r^{k_1}_{P^\prime P}
\end{align*}

Thus property (4) holds.
\end{proof}

Clearly, if $k,l\ge M(P,P^\prime)$ the above proposition shows that
$r=r^k_{P^\prime P}$ and $s=r^l_{PP^\prime}$ form a standard shift
equivalence between $f_P$ and $f_{P^\prime}$ of lag $k+l$.

Suppose $S$ is an isolated invariant set.  By Theorem \ref{equivthm}
we know that for any filtration pair $P=(N,L)$ the shift equivalence
class of $f_{P}:N_{L}\to N_{L}$ is an invariant for $S$.  This
invariant is, however, too restrictive for our desires.  Instead
consider the homotopy class of base point preserving maps on $N_{L}$
with $f_{P}$ as a representative.  We denote this collection
$h_{P}(S)$.

Suppose $X$ and $Y$ are pointed topological spaces and $[f]:X\to X$
and $[g]:Y\to Y$ are homotopy classes of base-point preserving maps.
Then the definition of shift equivalence states that $(X,[f])$ and
$(Y,[g])$ are shift equivalent if there exist homotopy classes of maps
$[r]:X\to Y$ and $[s]:Y\to X$ such that $[g\circ r]=[r\circ f]$,
$[s\circ g]=[f\circ s]$, $[r\circ s]=[g^{m}]$ and $[s\circ r]=[f^{m}]$
for some $m$.

We may now make the following definition.

\begin{defn}\label{homotopy_index}
Let $S$ be an isolated invariant set for a continuous map $f$.  Then
define the {\em discrete homotopy Conley index} of $S$, $h(S)$, to be the shift
equivalence class of $h_{P}(S)$ where $P=(N,L)$ is a filtration pair
for $S$.
\end{defn}

Once this invariant has been defined, we may apply functors to it to
obtain new invariants.  For instance, we may take homology to obtain
the {\em homological Conley index}, $\Con_{*}(S)$.  In other words,
$\Con_{*}(S)$ is the shift equivalence class of $(f_{P})_{*}:
H_{*}(N_{L},[L])\to H_{*}(N_{L},[L])$ where $P=(N,L)$ is a filtration
pair for $S$.  Similarly, we define the {\em cohomological Conley
index}, $\Con^{*}(S)$.

For some situations we may want the Conley index to be a graded
abelian group together with a distinguished automorphism.  Thus, we
may apply a functor to the shift equivalence class $\Con_*(S)$ or
$\Con^*(S)$ such as the Leray functor.  One should consult \cite{Mr}
for details.

The Conley index has several important properties.

\begin{thm} [Continuation property]\label{continproperty}
Suppose $f_\lambda:X\times[0,1]\to X$ is a continuous homotopy and $N$
is an isolating neighborhood for $f_0$ with $S_0=\Inv (N,f_0)$.  Then
there is an $\ep>0$ such that $N$ is an isolating neighborhood for
every $f_\lambda$, $\lambda<\ep$.  Moreover,
$h(S_0,f_0)=h(S_\lambda,f_\lambda)$ where
$S_\lambda=\Inv(N,f_\lambda)$ and $\lambda<\ep$.
\end{thm}

\begin{proof}
This result is an easy consequence of Theorem~\ref{fpexist}.
\end{proof}

Moreover, if our space is actually a manifold we have the following
stronger version of Theorem \ref{continproperty}.

\begin{cor}
Let $f:M\to M$ be a continuous map of a manifold with isolating
neighborhood $N$ and isolated invariant set $S=\Inv N$.  There is a
neighborhood of $f$ in the $C_{0}$ topology with the property that $N$
is an isolating neighborhood for every $\tilde f$ in this
neighborhood.  Moreover, $h(S,f)=h(\tilde S,\tilde f)$ where $\tilde
S=\Inv(N,\tilde f)$.
\end{cor}

\begin{proof}
We use the fact that we may calculate $h(S,f)$ using a filtration pair
homeomorphic to a finite simplicial complex.  This implies that for
any $\tilde f$ sufficiently close to $f$ the maps $f$ and $\tilde f$
are homotopic by a homotopy through maps close to $f$.
\end{proof}

The corollary given above does not hold for arbitrary locally compact
metric spaces.  For instance, one may take
$X=\{\frac{1}{n}:n\in\Z^+\}\cup\{0\}$, $f=id$, $N=X$ and
$L=\emptyset$.  Then, for any $\ep>0$ define a map $f_\ep:X\to X$ by
putting $f_\ep(x)=x$ when $x\ge\ep$ and $f_\ep(x)=0$ otherwise.  It is
clear that the maps $f$ and $f_\ep$ are not shift equivalent.

We say that $h(S)=0$ if each pointed space map for $S$ is shift
equivalent (in the homotopy category) to the map on the pointed space
consisting of one point.

\begin{thm}
[Wa\.{z}ewski property] The set $\emptyset$ is an isolated invariant
set with $h(\emptyset)=0$.  In particular, the contrapositive states
that if $h(S)\ne 0$ then $S\ne\emptyset$.
\end{thm}
\begin{proof}
This result follows immediately {from} Definitions~\ref{filpair} and
\ref{homotopy_index}.
\end{proof}

\section{Attractors}

Recall that for any $x\in X$ we define the $\omega$-limit set of $x$
to be \[\omega(x)=\underset{N>0}{\bigcap}
\cl\Big(\underset{n>N}{\bigcup} \{f^n(x)\}\Big ).\]

We would like to do the same for backward iterates.  But, because
there is no unique backward trajectory of a point, we must modify the
definition; we must speak of the backward limit of a solution.  For
any solution $\sigma:\Z^{-}\to X$ define \[\alpha(\sigma)=
\underset{N>0}{\bigcap}\cl\Big(\underset{n>N}{\bigcup}\{\sigma
(-n)\}\Big)\]

A set $A\subset X$ is called an {\em attractor} if there exists a
compact neighborhood $U$ of $A$ such that \[f(U)\subset\Int U\] and
\[A=\bigcap_{n>0}f^{n}(U).\] The set $U$ is called an {\em attracting
neighborhood} of $A$.  If $X$ is invariant then define the {\em dual
repeller} $A^*$ to be
\[
A^*=\{x\in X:\omega(x)\cap A=\emptyset\}.
\]
If $X$ is not invariant then we only consider $x\in\Inv X$.  The pair
$(A,A^*)$ is called an {\em attractor-repeller decomposition} of $X$.
For any two sets $B,C\in X$ we define $C(B,C;X)$, the {\em connecting
orbits} {from} $B$ to $C$, to be
\[
\{x\in X\setminus(B\cup C): \exists\text{ solution }\sigma:\Z\to X\text{
with }\alpha(\sigma)\subset B, \omega(x)\subset C\}.
\]

\begin{prop}
If $(A,A^{*})$ is an attractor-repeller decomposition of a compact
invariant space $X$ then
\begin{enumerate}
\item
$X=A\cup A^{*}\cup C(A^{*},A;X)$,
\item
there exists $\ep>0$ such that there exists no $\ep$-chain {from} $A$
to $A^{*}$, and
\item
$A$ and $A^{*}$ are disjoint isolated invariant sets.
\end{enumerate}
\end{prop}
\begin{proof}
The first of these properties is immediate {from} the definitions.
The second follows {from} the fact that if $U$ is an attracting
neighborhood for $A$ and $\ep < d( f(U), \cl (X\setminus U))$ then
there exists no $\ep$-chain {from} a point of $f(U)$ to a point of
$\cl (X \setminus U)$.  The only non-trivial part of the third
property is that $A^{*}$ is an isolated invariant set.  This holds
because whenever $U$ is an attracting neighborhood for $A$, then $\cl
(X\setminus U)$ is an isolating neighborhood for $A^{*}$.
\end{proof}

\begin{lem}
\label{attnbdlem}
If $V$ and $W$ are attracting neighborhoods then $V\cap W$, $V\cup W$
and $f^{-1}(V)$ are attracting neighborhoods.
\end{lem}
\begin{proof}
These results are immediate except for the fact that $f^{-1}(V)$ is an
attracting neighborhood.  To see this observe that since $f(V) \subset
\Int V,$
\[
f(f^{-1}(V)) \subset V \subset f^{-1}(f(V)) \subset f^{-1}(\Int V)
\subset \Int f^{-1}(V).
\]
\end{proof}

The following lemma is immediate if $f$ is a homeomorphism, but it
requires proof otherwise.

\begin{lem}
\label{close_nbhd}
Suppose $A$ is an attractor and $V$ is a neighborhood of $A$ in 
some manifold.  Then there is an attracting neighborhood $U\subset V$.
\end{lem}
\begin{proof}
It suffices to prove the result for an attractor $A$ consisting of a
single point since we can always consider the map induced by $f$
on the quotient space $X/A$ which has the single point $[A]$ as
attractor.

Thus we assume $A$ consists of one point.  Let $N$ be an attracting
neighborhood for $A$, so $P = (N, L)$ with $L = \emptyset$ is a
filtration pair for the invariant set $A.$ The associated pointed
space $N_L = N \cup [L]$ consists of $N$ together with a disjoint
single point $[L].$ The map $f_P : N_L \to N_L$ has precisely two
fixed points, $A$ and $[L]$.

Given any neighborhood $V$ of $A$ it follows {from}
Theorem~\ref{small_fp} that there is filtration pair $Q = (Y,K)$ for
$A$ inside $V$.  Let $r:Y_K \to N_L$ be one of the semi- conjugacies
{from} the shift equivalence between $f_Q$ and $f_P.$ Then $U_0 =
r^{-1}(N_L \setminus \{[L]\})$ is an open and closed subset of $Y_K$ which is
forward invariant, contains $A$, and is disjoint {from} $K.$ It follows that
$f_Q( U_0) \subset \Int U_0$ so if $p: Y \to Y_K$ is the quotient map
then $U = p^{-1}(U_0)$ is an attracting neighborhood for $A$ which is
contained in $V$.
\end{proof}

We would like to speak of an attractor-repeller decomposition of an
isolated invariant set $S$.  Clearly, all the definitions and results
given above apply.  However, we would like to define our attracting
neighborhood as a subset of the ambient space $X$, not as a subset of
$S$.  An obvious definition of such an attracting neighborhood $U$
would be that $U$ be a neighborhood of $A$ in $X$ with $f(U\cap
S)\subset \Int U\cap S$ and $A=\bigcap_{n>0}f^{n}(U\cap S)$.  In fact,
as we shall see later, given a filtration pair $P=(N,L)$ for $S$, we
may choose $U\subset N$ so that $U$ is a neighborhood of $L$ and in
the pointed space $N_{L}$, $f_{P}(U)\subset\Int U$.

\section{Morse decompositions}

Suppose $\p$ is a finite set.  We say that a relation $<$ is a {\em
partial order} on $\p$ provided for all $p,q\in\p$
\begin{enumerate}
\item
$p<p$ never holds, and
\item
$p<q$ and $q<r$ implies $p<r$.
\item
If, in addition, for any $p\ne q$ we have $p<q$ or $q<p$ then $<$ is
called a {\em total ordering}.
\end{enumerate}

Suppose $\p$ has a partial order $<$.  A subset $I\subset\p$ is called
an {\em interval} if $p<r<q$ and $p,q\in I$ then $r\in I$.  An
interval $I$ is called {\em attracting} provided $r<q$ and $q\in I$
implies $r\in I$.  We will denote the attracting interval with unique
maximal element $p$ by $I_{p}$.  We denote the collection of intervals
and attracting intervals by $\I=\I(\p,<)$ and $\A=\A(\p,<)$
respectively.

\begin{defn}
\label{morse_decomp}
The collection of disjoint isolated invariant sets
$\M(\p,<)=\{M_{p}\subset S:p\in\p\}$ is a {\em Morse decomposition} of
an isolated invariant set $S$ if for every $x\in S$ and every solution
$\sigma:\Z\to S$ for $x$ we have either
\begin{enumerate}
\item
$\sigma(\Z)\subset M_{p}$ for some $p\in\p$, or 
\item
$\omega(\sigma)\subset M_{p}$ and $\alpha(\sigma)\subset M_{q}$ for
some $p<q$.
\end{enumerate}
For any interval $I\subset\p$ define the set
\[
M_{I}=\bigcup_{p\in I}M_{p}\cup\bigcup_{q,r\in I}C(M_{q},M_{r};S).
\]
\end{defn}

Notice that the partial ordering on the Morse sets is analogous to
Smale's no-cycle property for basic sets of an Axiom A diffeomorphism
\cite{S}.

\begin{prop}
\label{attr_nbhd}
Suppose $f(X) \subset \Int X$ so $S = \Inv X$ is an attractor in $X$.
Let $\M(\p,<)=\{M_{p}\subset S:p\in\p\}$ be a Morse decomposition of
$S$.  If $I\subset \p$ is an attracting interval then $M_{I}$ is an
attractor, i.e.  there is a compact neighborhood $U_I\subset X$ of
$M_I$ such that
\begin{enumerate}
\item
$f(U_I) \subset \Int U_I$, and
\item
$M_I = \Inv(U_I) = \bigcap_{n>0}f^n(U_I)$.
\end{enumerate}
Moreover, if $V_I$ is an open neighborhood of $M_I$ then we can find
$U_I\subset V_I$.
\end{prop}

\begin{proof}
An attracting interval $I$ in $\p$ for the order $<$ is also an
attracting interval for a total ordering of the set $\p$ which is
compatible with $<$.  Thus, without loss of generality we may assume
that $<$ is a total order.

We will give the proof by induction on the number of elements $n$ in
$\p$.  Clearly if $n = 1$ then $M_I = S$ so we may take $U_I = X.$
Hence we may assume the result for $n =r-1$ and need to prove it for
$n=r$.

Let $p_r$ be the maximal element of $\p$ and define $J$ to be the
attracting interval $\p \setminus \{p_r\}$.  Let $B = \bigcup_{p \in
J} M_p$ and choose a compact neighborhood $V$ of $B$ which is disjoint
{from} $M_{p_r}.$ We wish to show that for some smaller neighborhood
$W$ of $B$ we have the property that $\omega(x) \subset B$ for all $x
\in W.$ We assume this is not the case and show that leads to a
contradiction.

If no such $W$ exists then there is a convergent sequence $\{x_n\}$
such that $\lim x_n \in B$ and $\omega(x_n)$ is not a subset of $B.$
{From} this it follows that we can choose $y_n$ with the properties
\begin{enumerate}
\item
$y_n \in V$,
\item
$f^k(y_n) \notin V$ for $0 < k < n,$ and
\item
$y_n = f^{j_n}(x_n)$ for some $j_n > 0.$
\end{enumerate}

By choosing a subsequence we may assume that the sequence $\{y_n\}$
converges, say to $z \in V.$ We claim that $z \in S$.  This follows
because $f^k(z) \in X$ for all $k > 0$ and for any $m >0$ the sequence
$\{w_n = f^{j_n-m}(x_n)\}$ is well defined for sufficiently large $n$
and has a subsequence which converges to a point of $f^{-m}(z).$ Thus
we have shown $z \in \Inv^m X$ for all $m>0$ so by Lemma~\ref{Invm} we
have $z \in \Inv X = S.$

Also, for every $k > 0$ $f^k(z) \notin \Int V$ since $f^k(z) =
\lim_{n \to \infty} f^k(y_n).$ Thus $\omega(z)$ is disjoint from
$\Int V$ and the definition of Morse decomposition implies that
$\omega(z) \subset \M_{p_r}$.  Let $\sigma$ be a solution for $z$.
Then $\alpha(z) \subset \M_q$ for some $q \in \p$, so $q \ge p_r.$ But
this leads to a contradiction since $p_r$ is maximal, so $q > p_r$ is
impossible, but $q = p_r$ would imply that $z \in M_{p_r}$ which is
false.

Thus we have shown that there is a compact neighborhood $W \subset V$
of $B$ with the property that $\omega(x) \subset B$ for all $x \in W.$
To each $x\in W$ we may associate $n(x)$ defined to be the smallest
non-negative integer such that $f^{n(x)}(x) \in \Int W.$ The function
$n:W\to\Z^{+}$ is upper semicontinuous.  Since $W$ is compact
$N=\max\{n(x):x\in W\}$ is finite.  Let $U_J=f^{-N}(W)$.  If $x \in
U_J$ then $f^k(x) \in \Int W$ for some $0 < k \le N$ so $f(x) \in
f^{-N}(\Int W) \subset \Int U_J.$ Hence $f(U_J) \subset \Int U_J.$
Since $M_{p_r} \cap U_J = \emptyset$ and $M_J \subset U_J$ it follows
that $M_J = \Inv U_J.$ The fact that $\Inv U_J = \cap_{n>0}f^n(U_J)$
follows {from} Lemma \ref{Invm}.

If $I=J$ then we're done.  If $I \ne J$ then $I$ is a subinterval of
$J.$ Note that $\M(J,<)$ is a Morse decomposition for the map $f: U_J
\to U_J$ with attractor $M_J$.  Since $J$ has one fewer elements than
$\p$ our induction hypothesis implies there is an attracting
neighborhood $U_I \subset U_J$ for $M_I.$

Finally, we remark that by Lemma~\ref{close_nbhd} we may assume that
$U_{I}\subset V_{I}$.
\end{proof}

Suppose $\M(\p,<)=\{M_{p}\subset S:p\in\p\}$ is a Morse decomposition
of an isolated invariant set $S$, and $P =(N,L)$ is a filtration pair
with associated pointed space map $f_P:N_L \to N_L$.  Then $N_L$ is an
attracting neighborhood for an attractor $\hat S = \Inv N_L =
\bigcap_{n>0}f^n_P(N_L).$ Moreover there is a Morse decomposition of
$\hat S$ associated to $\M(\p,<)$ obtained by adding the single point
$[L]$ as a minimal Morse set.  More precisely, we define the
associated Morse decomposition as follows.

\begin{defn}
Let $\hat \p = \p \cup \{0\}$ and specify $0 < p$ for all $p \in \p.$
Define the {\em associated Morse decomposition} to $\M(\p,<)$ to be
\[
\hat \M(\p,<)=\{\hat M_{p}: \hat M_p = M_p, p\in \p \} \cup \{M_0\},
\]
where $M_0=\{[L]\}.$
\end{defn}

It is not immediately clear that the associated Morse decomposition is
in fact a Morse decomposition, so we prove that now.

\begin{prop}
\label {assoc_md} Suppose $\M(\p,<)$ is a Morse decomposition for the
isolated invariant set $S$ with filtration pair $(N,L)$.  Then the
associated Morse decomposition $\hat \M(\hat \p,<)$ is a Morse
decomposition for the attractor $\hat S = \Inv N_L.$
\end{prop}

\begin{proof}
The sets $\hat M_p$ with $p \in \hat \p$ can immediately be seen to
satisfy the properties of the definition of Morse decomposition
\ref{morse_decomp} with one exception.  If $z \in \hat S$ with
$\omega(z) = \{[L]\}$ and $\sigma$ is a solution for $z$ then it is
not immediate that $\alpha(\sigma) \subset \hat M_q$ for some $q \in
\hat \p.$ However, the proof of Proposition \ref{attr_nbhd} did not
make use of this property so the conclusion of this proposition is
valid for $\hat \M(\hat \p,<)$.  In particular if $I\subset\p$ is an
attracting interval then $M_I$ is an attractor.

We now consider $z \in \hat S$ with $\omega(z) = \{[L]\}$ and let
$\sigma$ be a solution for $z$.  Since $\alpha(\sigma) \subset S$ and
$S$ has a Morse decomposition given by $\{\hat M_p : p \ne 0\}$ it
must be the case that either $\alpha(\sigma)$ is a subset of a single
$\hat M_q$ or $\alpha(\sigma)$ contains points of $\hat M_q$ and $\hat
M_p$ with $p \ne q.$ But if the attracting interval $I$ is chosen to
contain one of $p$ and $q$ and not the other then the existence of the
attracting neighborhood $U_I$ makes it impossible for $\alpha(\sigma)$
to contain points of both $\hat M_q$ and $\hat M_p$.  Thus
$\alpha(\sigma)\subset \hat M_{q}$ for some $q\in\p$.
\end{proof}

\begin{cor}
If $\M(\p,<)$ is a Morse decomposition of an isolated invariant set
$S$ then for any interval $I\subset\p$ the set $M_{I}$ is an isolated
invariant set.
\end{cor}
\begin{proof}
Without loss of generality we may assume that the hypotheses of
Proposition~\ref{attr_nbhd} hold, for if not, then by
Proposition~\ref{assoc_md} we may simply consider the pointed space
for $S$ and the associated Morse decomposition there.  Moreover, since
any interval is also an interval for a total ordering we may assume
$(\p, <)$ is a total order.  In this order there are attracting
intervals $J_0$ and $J_1$ such that $J_1 = J_0 \cup I.$ Consider the
attracting neighborhood $U_1$ for $\M_{J_1}$ given by Proposition
\ref{attr_nbhd} and the map $f: U_1 \to U_1.$ A Morse decomposition is
given by $\M(J_1,<)$.  Another application of Proposition
\ref{attr_nbhd} implies that $\M_{J_0}$ is an attractor with
attracting neighborhood $U_0 \subset U_1.$ Clearly, with respect to
$f: U_1 \to U_1$ the set $M_I$ is the dual repeller to $M_{J_0}$ and
hence is closed.  The set $N = U_1 \setminus \Int U_0$ is an isolating
neighborhood of $M_I.$
\end{proof}

\section{Morse set filtrations}

A variation of the following definition was introduced by Franzosa to
study Morse decompositions for continuous time dynamical systems
\cite{F}.  It is also similar to the neighborhood networks introduced
by Robbin and Salamon \cite{RS} for discrete dynamical systems.

\begin{defn} A collection of compact sets $\N(\p,<)=\{N(I)\subset
X:I\in\A\}$ is called an {\em Morse set filtration} for $\M(\p,<)$
provided that for any attracting intervals $I,J\subset\p$ the
following conditions hold.
\begin{enumerate}
\item
$(N(I),N(\emptyset))$ is a filtration pair for $M_{I}$,
\item
$N(I)\cap N(J)=N(I\cap J)$, and
\item
$N(I)\cup N(J)=N(I\cup J)$.
\end{enumerate}
\end{defn}

\begin{prop}
Suppose $\N(\p,<)$ is a Morse set filtration for a Morse decomposition
$\M(\p,<)$ and $J\subset I$ are attracting intervals.  Then
$(N(I),N(J))$ is a filtration pair for $M_{I\setminus J}$.
\end{prop}
\begin{proof} {From} property (2) of the definition above it is clear
that $N(I) \subset N(J)$ whenever $I$ and $J$ are attracting intervals
with $J\subset I$.  It is then straightforward to check that
$(N(I),N(J))$ satisfies the requirements of Definition~\ref{filpair}.
\end{proof}

\begin{lem}
\label{attnbhds}
Suppose $f(X) \subset \Int X$ and $S = \Inv X$ has a Morse
decomposition $\M(\p,<)=\{M_{p}\subset S:p\in\p\}$.  Then there exists
a collection \[\{V(I): V(I)\text{ is an attracting neighborhood for
}M_{I},\ I\in\A\}\] with the property that $V(I)\cap V(J)=V(I\cap J)$
and $V(I)\cup V(J)=V(I\cup J)$.
\end{lem}
\begin{proof}
We shall build the $V(I)$ inductively.  The induction hypothesis is as
follows.  Suppose $J\subset \p$ is an attracting interval and that for
each attracting interval $I\subset J$ there is an attracting
neighborhood $V(I)$ satisfying the conclusions of the lemma.
Moreover, if $p\in J$ and $q\not\in J$ with $p\not<q$ then
$\cl(V(I_{p})\setminus V(I_{p}\cap I_{q}))\cap M_{I_{q}}=\emptyset$.

For $J=\emptyset$ take $V(\emptyset)=\emptyset$.  Suppose the
induction holds for an attracting interval $J$ and $q$ is minimal in
$\p\setminus J$.  We will define $V(I)$ for all $I$ in the attracting
interval $J\cup\{q\}$.  Without loss of generality we may assume that
for each $p\ne q$ we have $M_{I_{p}}\cap M_{I_{q}}\subset V(I_{p}\cap
I_{q})$.  For if not then we may replace each $V(I)$ by $f^{-N}(V(I))$
for some $N$ sufficiently large.

By the induction hypothesis and the above remark, there exists an open
neighborhood $U$ of $M_{I_{q}}$ such that $U$ is disjoint {from}
$\cl(V(I_{p})\setminus V(I_{p}\cap I_{q}))$ for all $p\in J$, $p\not<q$ and
that $U$ is disjoint {from} $\cl(M_{I_{p}}\setminus V(I_{p}\cap I_{q}))$ for
all $p\not < q$.  By Lemma~\ref{close_nbhd} there is an attracting
neighborhood $V\subset U$ for $M_{I_{q}}$.  Define
\[V(I_{q})=V\cup(\bigcup_{p\in I_{q}\setminus\{q\}}V(I_{p}))\] and more
generally for any attracting interval $I\subset J\cup \{q\}$
\[V(I)=\bigcup_{p\in I}V(I_{p}).\] It is straightforward to check that
this construction completes the inductive step.
\end{proof}

\begin{thm}
Suppose $\M(\p,>)$ is a Morse decomposition of an isolated invariant
set $S$.  Then there exists a Morse set filtration for $\M(\p,<)$.
\end{thm}
\begin{proof}
Let $P=(N,L)$ be a filtration pair for $S$ and $\hat S=\Inv N_{L}$.
Let $\hat\M(\hat\p,<)$ denote the associated Morse decomposition of
$\hat S$.  By Lemma \ref{attnbhds} there exists a collection
\[\{V(I)\subset N_{L}: V(I)\text{ is an attracting neighborhood for
}M_{I},\ I\in\A(<)\}\] with the property that $V(I)\cap V(J)=V(I\cap
J)$ and $V(I)\cup V(J)=V(I\cup J)$.  To obtain a Morse set filtration
for $\M(\p,<)$ we simply take $N(I)=q^{-1}(V(I\cup\{0\}))$ where
$0\in\hat\p$ corresponds to the Morse set $\{[L]\}$ and $q:N\to N_{L}$
is the quotient map.  It is clear that this collection forms a Morse
set filtration.
\end{proof}

\section{Shift equivalence and Szymczak's category}

In \cite{Sz} Szymczak provides a different, categorical, definition of
the Conley index.  We conclude this discussion by showing that his
definition is identical to ours.

For any category $\K$ define the category of {\em objects equipped
with a morphism over} $\K$, denoted $\K_{m}$, as follows.

\[
\Ob(\K_{m})=\{(X,f):X\in\Ob(\K), f\in\Mor_{\K}(X,X)\},
\]
\[
\Mor_{\K_{m}}((X,f),(X^{\prime},f^{\prime}))= M((X,f
),(X^{\prime},f^{\prime}))/\sim
\]
where
\[
M((X,f),(X^{\prime},f^{\prime}))= \{
(g,n)\in\Mor_{\K}(X,X^{\prime})\times\Z^{+}: g f=f^{\prime} g\}
\]
and $\sim$ is the following equivalence relation.  $(g_{1},n_{1})\sim
(g_{2},n_{2})$ if and only if there exists $k\in\Z^{+}$ such that the
following diagram commutes.
\begin{equation*}
\begin{CD}
X @>f^{n_{1}+k}>> X\\ @Vf^{n_{2}+k}VV @VVg_{2}V\\ X @>g_{1}>>
X^{\prime}.
\end{CD}
\end{equation*}
Henceforth we shall call this category the {\em Szymczak category}.

\begin{prop}
\label{shiftprop}
Suppose that $(X,f),(X^{\prime},f^{\prime})\in\Ob(\K_{m})$. Then
$(X,f)$ and $(X^{\prime},f^{\prime})$ are isomorphic in the Szymczak
category if and only if they are shift equivalent.
\end{prop}
\begin{proof}
Suppose $f$ and $f^{\prime}$ are shift equivalent with $rs=f^{m}$ and
$sr=(f^{\prime})^{m}$.  We claim that
$[r,m]\in\Mor_{\K_{m}}((X^{\prime},f^{\prime}),(X,f))$ is an
isomorphism with inverse
$[s,0]\in\Mor_{\K_{m}}((X,f),(X^{\prime},f^{\prime}))$.  Because the
diagram
\begin{equation*}
\begin{CD}
X @>f^{m}>> X\\ @VidVV @VVidV\\ X @>s r>> X.
\end{CD}
\end{equation*}
commutes we see that $(s,0)\circ(r,m)=(s r,m)\sim(id_{X},0)$.
Similarly, $(rs,m)\sim(id_{X^{\prime}},0)$.  Thus, $[r,m]$ and $[s,0]$
are inverse isomorphisms.

Now suppose that $(X,f)$ and $(X^{\prime},f^{\prime})$ are isomorphic
objects in Szymczak's category.  In particular, suppose
$[r,t]\in\Mor_{\K_{m}}((X,f),(X^{\prime},f^{\prime}))$ and
$[s,u]\in\Mor_{\K_{m}}((X^{\prime},f^{\prime}),(X,f))$ are inverse
isomorphisms.  We have $(s r, t+u)\sim(id_{X},0)$ and $(r
s,t+u)\sim(id_{X^{\prime}},0)$.  Thus, we have $r f=f^{\prime} r$, $f
s=s f^{\prime}$ and $k_{1},k_{2}\in\Z^{+}$ such that the following
diagrams commute.
\begin{multicols}{2}
\begin{center}
$\begin{CD} X @>f^{t+u+k_{1}}>> X\\ @Vf^{k_{1}}VV @VVidV\\ X @>s r>> X
\end{CD}$\\

$\begin{CD} X^{\prime} @>(f^{\prime})^{t+u+k_{2}}>> X^{\prime}\\
@V(f^{\prime})^{k_{2}}VV @VVidV\\ X^{\prime} @>r s>> X^{\prime}.
\end{CD}$
\end{center}
\end{multicols}
Let $k=\max\{k_{1},k_{2}\}$, $m=t+u+k$, $s^{\prime}=s$ and
$r^{\prime}=r f^{k}$.  So defined, we have $f s^{\prime}=s^{\prime}
f^{\prime}$, $r^{\prime} f=f^{\prime} r^{\prime}$, $r^{\prime}
s^{\prime}=(f^{\prime})^{m}$ and $s^{\prime} r^{\prime}=f^{m}$.  Thus
$f$ and $f^{\prime}$ are shift equivalent.
\end{proof}

Let $Htop_{*}$ be the category of pointed topological spaces with
homotopy classes of base-point preserving maps.  An immediate
consequence of Proposition \ref{shiftprop} is that Szymczak's
definition of the Conley index is equivalent to ours.  We state this
explicitly in the following proposition.

\begin{prop}
Suppose $P=(N,L)$ is a filtration pair for an isolated invariant set
$S$.  Then the isomorphism class of $(N_{L},f_{P})$ in the Szymczak
category $(Htop_{*})_{m}$ is precisely $h(S)$, the Conley index of
$S$.
\end{prop}

\end{document}